\documentclass[hidelinks]{siamart190516}
\usepackage{graphicx,amsfonts,upquote,verbatim}
\def\pput(#1,#2)#3{\noindent\smash{\raise#2pt\hbox to 0pt
   {\kern #1pt #3\hss}}\ignorespaces}

\def\Re{\hbox{\rm Re\kern .3pt}}

\def\pexact{p_{\hbox{\footnotesize\rm exact}}^{}}
\title{AAA interpolation of equispaced data}
\author{Daan Huybrechs\thanks{\texttt{daan.huybrechs@kuleuven.be},
Dept.~Computer Sci., KU Leuven, Celestijnenlaan 200A, BE-3001 Leuven, Belgium}
\and
Lloyd N.~Trefethen\thanks{\texttt{trefethen@maths.ox.ac.uk},
Mathematical Institute, University of Oxford, Oxford OX2 6GG, UK.}}

\begin{document}

\maketitle


\begin{abstract}
We propose AAA rational approximation as a method for 
interpolating or approximating smooth functions from 
equispaced data samples.
Although it is always better to approximate from
large numbers of samples if they are available, whether equispaced or not,
this method often performs impressively even when the sampling grid is
fairly coarse.
In most cases it gives more accurate approximations than other
methods.
\end{abstract}


\pagestyle{myheadings}
\thispagestyle{plain}
\markboth{\sc Huybrechs and Trefethen}
{\sc AAA approximation of equispaced data}

\section{Introduction}

The aim of this paper is to propose a method for interpolation
of real or complex data in equispaced points on an interval, which without
loss of generality we take to be $[-1,1]$.  In its basic form the method
simply computes a AAA rational 
approximation\footnote{pronounced ``triple-A''}~\cite{aaa}
to the data, and thus the interpolant is a numerical one,
not mathematically exact: a crucial advantage for
robustness.  In Chebfun~\cite{chebfun}, the fit can be computed
to the default relative accuracy $10^{-13}$ by the command
\begin{equation}
{\tt r = aaa(F)},
\label{command1}
\end{equation}
where $F$ is the vector of data values.
(As explained in Section~3, an adjustment is made if
the AAA approximant turns out to have
poles in the interval of approximation.)
If interpolation by a polynomial rather than a rational function is desired,
this can by determined by a further step in which $r$ is approximated by a Chebyshev series,
\begin{equation}
{\tt p = chebfun(r)}.
\label{command2}
\end{equation}

\begin{figure}
\begin{center}
\vspace*{15pt}
\includegraphics[scale=.9]{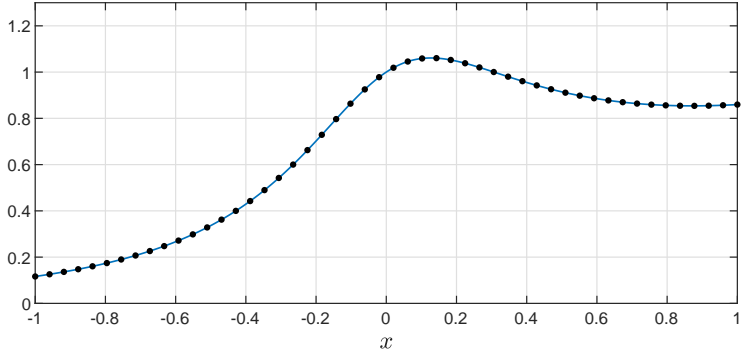}
\vspace*{5pt}
\end{center}
\caption{\small\label{fig1}AAA interpolation of
$f(x) = e^x/\sqrt{1+9x^2}$ in $50$ equispaced points in $[-1,1]$.
The rational interpolant, of degree $17$, matches $f$ to accuracy $3.3 \times 10^{-14}$
at the sample points and $9.6\times 10^{-14}$ on $[-1,1]$.}
\end{figure}

For example, Figure~\ref{fig1} shows the AAA interpolant $r$ of
$f(x) = e^x/\sqrt{1+9x^2}$ in 50 equispaced points of $[-1,1]$.
This is a rational function of degree $17$ with accuracy $\|f-r\|
\approx 9.6\times 10^{-14}$, computed in about a millisecond on 
a laptop.  (Throughout this paper, $\|\cdot\|$ is the
$\infty$-norm over $[-1,1]$.)  The Chebfun polynomial approximation $p$
to $r$ has degree $104$ and the same accuracy $\|f-p\kern .8pt \| \approx
9.6\times 10^{-14}$.  The exact degree 49 polynomial interpolant
$\pexact$ to the data, by contrast, has error $\|f-\pexact\| \approx
109.3$ because of the Runge phenomenon~\cite{runge,ATAP}.
It is fascinating that one can generate high-degree polynomial interpolants
like this that are so much better between the sample points 
than polynomial interpolants of minimal degree, but we see no
particular advantage in $p(x)$ as compared with $r(x)$, so for the remainder
of the paper, we just discuss $r(x)$.

The problem of interpolating or approximating equispaced data arises
in countless applications, and there is a large literature on the
subject, with many algorithms having been proposed, some of them
highly effective in practice.  One reason why no single algorithm
has taken over is that there is an unavoidable tradeoff in this
problem between accuracy and stability.  In particular, if $n$
equispaced samples are taken of a function $f$ that is analytic
on $[-1,1]$, one might expect that exponential convergence to $f$
should be possible as $n\to\infty$.  However, the {\em impossibility
theorem\/} asserts that exponential convergence is only possible in
tandem with exponential instability, and that, conversely,
a stable algorithm can converge at best at a root-exponential rate
$\|f-r_n\| = \exp(-C\sqrt n\kern 1pt)$, $C>0$~\cite{ptk}.  In practice, it is
usual to operate in an in-between regime, accepting some instability
as the price of better accuracy.  In the face of this complexity, it
follows that different algorithms may be advantageous for different
classes of functions, and that pinning down the properties of any
particular algorithm may not be straightforward.

In this complicated situation we will do our best to elucidate the
properties of AAA interpolation.  First we compare its
performance numerically against that of some existing algorithms for a collection
of approximation problems (sections~2 and~3), and then we present certain
theoretical considerations (section 4).  The final section briefly
discusses AAA variants, the effect of noise, and other issues.

Apart from some remarks around Figure~\ref{twoposs},
we will not describe details of the AAA algorithm, partly because
these can be found elsewhere~\cite{aaa} and partly because the
essential point here is not AAA per se but just 
rational approximation. 
At present, AAA appears to be the best general-purpose
rational approximation tool available, but other methods may come along
in the future.

\section{Existing methods}

Many methods have been proposed for interpolation or approximation of equispaced
data, and we will not attempt a comprehensive review.  We will,
however, mention the main categories of methods and choose five specific
examples for numerical comparisons.  For previous surveys with further
references, see~\cite{bo1,bo2,ptk}.

With any interpolant or approximant, there are always the two
questions of the {\em form\/} of the approximant and the {\em
method\/} of defining or computing it.  The main forms that have been
advocated are polynomials or piecewise polynomials,
Fourier series, rational functions, radial
basis functions (RBFs), and various modifications and combinations
of these.  The methods proposed generally involve mathematically
exact interpolation or some version of least-squares approximation.
(Ironically, because of conditioning issues, exact interpolants may
be less accurate in floating-point arithmetic than least-squares
approximations, even at the sample points, let alone in-between.)
Almost every method involves a choice of one or two parameters, which usually
affect the tradeoff between accuracy and stability and can
therefore be interpreted in part as regularization parameters.
AAA interpolation may appear at first as an exception to this rule,
but a parameter implicitly involved is the tolerance, which in
the Chebfun implementation is set by default to $10^{-13}$.  We will
discuss this further in Section 4.

\smallskip
{\em Polynomial least-squares.}
Interpolation of $n$ data values by a polynomial of degree $n-1$ leads to
exponential instability at a rate $O(2^n)$, as has been known
since Runge in 1901~\cite{runge,ATAP}.  Least-squares fitting
by a polynomial of degree $d< n-1$, however, is better behaved.
To cut off the exponential growth as $n\to\infty$ entirely, $d$ must
be restricted to size $O(\sqrt n\kern 1pt )$~\cite{bx,rakh}, but one
can often get away with larger values in practice, and a simple
choice is $d  \approx n/\gamma$, where $\gamma>1$ is an {\em oversampling ratio}.
According to equation (4.1) of~\cite{ptk}, this cuts the exponential unstable
growth rate from $2^n$ to $C^n$ with
\begin{equation}
C = \left[(1+\alpha)^{1+\alpha} (1-\alpha)^{1-\alpha}\right]^{1/2}, \quad
\alpha =  1/\gamma.
\label{kuijlaars}
\end{equation}
For example, $\gamma=2$ yields a growth rate of about $(3^{3/4}/2)^n \approx
(1.14)^n$, which is mild
enough to be a good choice in many applications.
Our experiments of Figure~\ref{figbig} in the
next section use this value $\gamma=2$.

\smallskip
{\em Fourier, polynomial, and RBF extensions.}
The idea of {\em Fourier extension\/} is to approximate $f$ by a Fourier series
tied not to $[-1,1]$ but to a larger domain
$[-T,T\kern 1pt]$~\cite{boyd02,bhp}.
The fit is carried out by least-squares or regularized least-squares,
often simply
by means of the backslash operator in MATLAB, as has been analyzed
by Adcock, Huybrechs, and Vaquero~\cite{ahmv,h09} and Lyon~\cite{lyon}.
A related idea is {\em polynomial extension}, in which $f$ is approximated by
polynomials expressed in a basis of orthogonal polynomials defined on
an interval $[-T,T\kern 1pt]$~\cite{as}.  A third possibility is {\em
RBF extension,} in which $f$ is approximated by smooth RBFs whose centers
extend outside $[-1,1]$~\cite{flt,piret}.  In Figure~\ref{figbig} of the
next section, we use Fourier extension with $T=2$ and an oversampling
ratio of $2$, so that the least-squares matrices have twice as many
rows as columns.

\smallskip
{\em Fourier series with corrections.}
If $f$ is periodic, trigonometric (Fourier) interpolation provides a
perfect approximation method: exponentially convergent and stable.
In the case of quadrature, this becomes the exponentially
convergent trapezoidal rule~\cite{tw}.  For nonperiodic $f$,
an attractive idea is to 
employ a trigonometric fit modified by corrections of one sort or another, 
often the addition of a polynomial term,
designed to mitigate the effect of the implicit discontinuity at
the boundary.  This idea goes back as far as James Gregory in 1670, before
Fourier analysis and even calculus~\cite{forn21}!\ \ The result
will not be exponentially convergent, but it can have an algebraic
convergence rate of arbitrary order depending on the choice of the
corrections, and the rate may improve to super-algebraic if the
correction order is taken to increase with $n$.  This idea has
been applied in many variations, an early example being a method
of Eckhoff with precursors he attributes to Krylov and Lanczos~\cite{eckhoff}.
In the ``Gregory interpolant''
of~\cite{jt}, an interpolant in the form of a sum of a trigonometric term and
a polynomial is constructed whose integral equals
the result for the Gregory quadrature formula.
Fornberg has proposed (for quadrature, not yet approximation)
a method of regularized endpoint corrections in
which extra parameters are introduced whose 
amplitudes are then limited by optimization~\cite{forn21}.
Figure~\ref{figbig} of the next section shows curves for
a least-squares method in which a Fourier series is combined with a polynomial term of
degree about $\sqrt n$, with an oversampling ratio of about $2$.

\smallskip
{\em Multi-domain methods.}
Related in spirit to methods involving boundary corrections are
methods in which $f$ is approximated by different functions over different
subintervals of $[-1,1]$---in the simplest case, a big central interval and two 
smaller intervals near the ends.  For examples see~\cite{boyd07,bo2,klein13,pg}.

\smallskip
{\em Splines.}  Splines, which are piecewise polynomials
satisfying certain continuity conditions, take the multi-domain idea further and
are an obvious candidate for
approximations that will not suffer from Gibbs oscillations at the boundaries.  
The most familiar case is cubic spline interpolants, where the sample
points are nodes separating cubic pieces with
continuity of function values and first and second derivatives.
Cubic splines (with the standard natural boundary conditions at the ends
and not-a-knot conditions one node in from the ends) are one of the methods
presented in Figure~\ref{figbig} of the next section.

\smallskip
{\em Mapping.}  By a conformal map, polynomial approximations can be
transformed to other approximations that are more suitable for equispaced
interpolation and approximation.  The prototypical method in this area
was introduced by Kosloff and Tal-Ezer~\cite{kte}, and there is also
a connection with prolate spheroidal wave functions~\cite{prolate}.
The general conformal mapping point of
view was put forward in~\cite{ht} and in~\cite[chapter 22]{ATAP}.
See also~\cite{boyd16}.

\smallskip
{\em Gegenbauer reconstruction.}  Another class of methods has been
developed from the point of view of edge detection and elimination of
the Gibbs phenomenon in harmonic analysis.  For entries
into this extensive literature, see~\cite{gt} and~\cite{tadmor}.

\smallskip
{\em Explicit regularization methods.}  Several other methods, often
nonlinear, have been proposed involving various strategies of
explicit regularization to counter the instability of
high-accuracy approximation~\cite{berzins,boyd92,chand,wmi}.
We emphasize that even many of the simpler numerical methods implicitly involve
regularization introduced by rounding errors, as will be discussed
in Section~4.

\smallskip
{\em Floater--Hormann rational interpolation:~Chebfun \verb|'equi'|.}
Finally, here is another method involving rational functions.
Floater and Hormann introduced a family of degree $n-1$ rational interpolants in
barycentric form whose weights can be adjusted
to achieve any prescribed order of accuracy~\cite{fh}.  (The AAA method also uses
a barycentric representation, but it is an approximant, in principle, not
an interpolant, and it is not very closely related to Floater--Hormann approximation.)
The method we show in Figure~\ref{figbig} of the next section, due to
Klein~\cite{gk,klein}, is based on interpolants
whose order of accuracy is adaptively determined via the
\verb|'equi'| option of the Chebfun constructor~\cite{bos,chebex}.

\section{Numerical comparison}

As mentioned in the opening paragraph, our
interpolation method consists of AAA approximation with its
standard tolerance of $10^{-13}$, so long as the approximant that
is produced has no ``bad poles,'' that is,
poles in the interval $[-1,1]$.  The principal
drawback of AAA approximation is that such poles sometimes appear---often
with such small residues that they do not contribute
meaningfully to the quality of the approximation (in which case
they may be called ``spurious poles'' or ``Froissart doublets''~\cite{ATAP}).
When the original
AAA paper~\cite{aaa} was published, a ``cleanup'' procedure was
proposed to address this problem.  We are no longer confident that
this procedure is very helpful, and instead, we now propose the
method of {\em AAA-least squares} (AAA-LS) introduced
in~\cite{costa}.  Here, if there are any bad poles, these are discarded, and
the other poles are retained to form the basis of a linear least-squares
fit to find a new rational approximation represented in partial fractions form.
For details, see the ``{\tt if any}'' block in the AAA part of
the code listed in the appendix.
Typically this correction makes little difference to accuracy down to 
levels of $10^{-7}$ or so, but it is often unsuccessful at achieving tighter
accuracies than this.

Poles in $[-1,1]$ almost never appear in the approximation of functions 
$f(x)$ that are complex (a case not illustrated here as it is less common
in applications).  For real problems, accordingly, another way of avoiding bad
poles is to perturb the data by a small, smooth complex function.
Overall, however, it must be said that the appearance
of unwanted poles in AAA approximants is not yet fully understood, and
it seems likely that improvements are in store in this active research area.

\begin{figure}[t]
\vspace*{5pt}
\begin{center}
\includegraphics[scale=1.1]{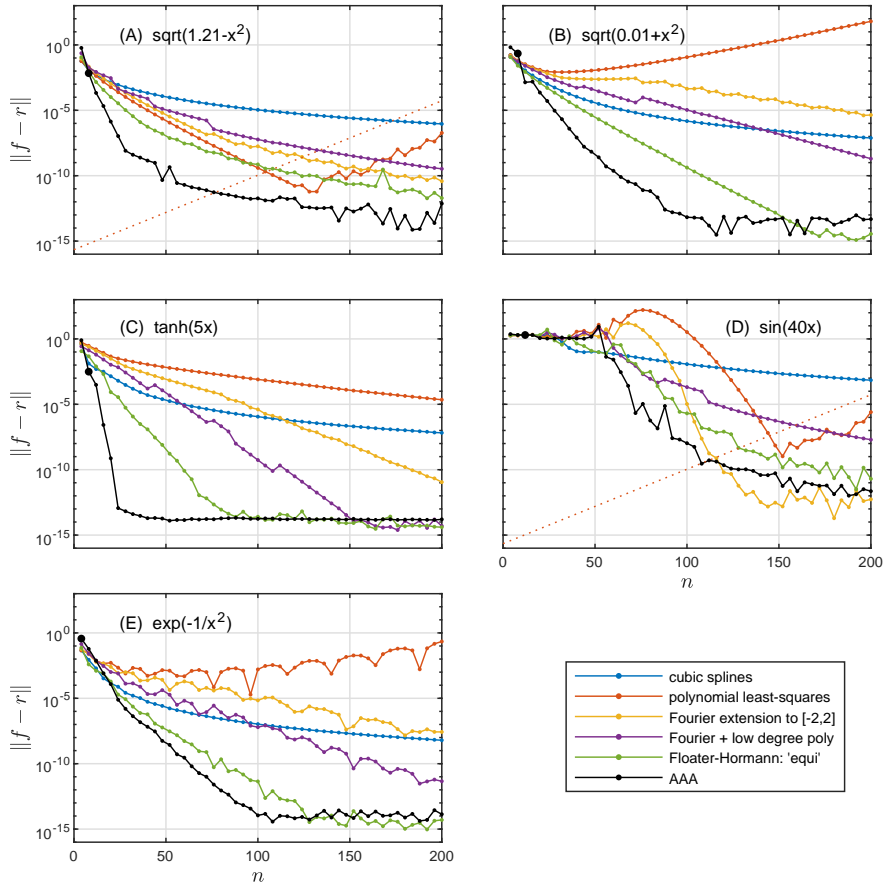}
\vspace{-3pt}
\end{center}
\caption{\small\label{figbig}Six approximation methods applied to five
smooth functions on $[-1,1]$.  In each case the horizontal axis is $n$,
the number of sample points, and the vertical axis is $\|f-r\|$, the maximum
error over $[-1,1]$.  The thicker dot on the AAA curve marks the final value of $n$
at which the rational function is an interpolant rather than just an approximation.
The dashed lines in (A) and (D) mark the instability estimate
$(1.14)^n\varepsilon_{\hbox{\scriptsize\rm machine}}$ from $(\ref{kuijlaars})$
with oversampling ratio $\gamma=2$.
The results are discussed in the text, and the code
is listed in the appendix.}
\end{figure}

Comparing the AAA method against other methods can quickly grow very
complicated since most methods have
adjustable parameters and there are any number of functions
one could apply them to.  To keep
the discussion under control, the panels of Figure~\ref{figbig}
correspond to five functions:
\begin{alignat}{3}
&f_A^{}(x) &&= \sqrt{1.21-x^2}\kern 20pt &&
 \hbox{(branch points at $\pm 1.1$),}\\[2pt]
&f_B^{}(x) &&= \sqrt{0.01+x^2}  && \hbox{(branch points at $\pm 0.1i$),}\\[3pt]
&f_C^{}(x) &&= \tanh(5x)        && \hbox{(poles on the imaginary axis),}\\[3.5pt]
&f_D^{}(x) &&= \sin(40x)        && \hbox{(entire, oscillatory),}\\[3pt]
&f_E^{}(x) &&= \exp(-1/x^2)     && \hbox{($C^\infty$ but not analytic).}
\end{alignat}
Each panel displays convergence curves for six methods:
\vskip 6pt
~~~~~~Cubic splines,
\vskip 3pt
~~~~~~Polynomial least-squares with oversampling ratio $\gamma = 2$,
\vskip 3pt
~~~~~~Fourier extension on $[-2,2]$ with oversampling ratio $\gamma = 2$,
\vskip 3pt
~~~~~~Fourier series plus polynomial of degree $\sqrt n$ with oversampling ratio $\gamma = 2$,
\vskip 3pt
~~~~~~Floater--Hormann rational interpolation: Chebfun \verb|'equi'|,
\vskip 3pt
~~~~~~AAA with tolerance $10^{-13}$.
\vskip 6pt
\noindent For details of the methods, see the code
listing in the appendix.  

Many observations can be drawn from Figure~\ref{figbig}.  The most
basic is that AAA consistently appears to be the best of the methods,
and is the first to reach accuracy $10^{-10}$ in every case.  It is
typical for AAA to converge twice as fast as the other methods, and for the
test function $f_C^{}(x) = \tanh(z)$, whose singularities consist of
poles that the AAA approximants readily capture, its superiority
is especially striking.

It is worth spelling out the meaning of the AAA convergence curves
of Figure~\ref{figbig}.
Each point on one of these curves corresponds to
a rational approximation whose error is $10^{-13}$ or less on the
discrete grid (at least if the partial fractions least-squares
procedure has not been invoked because of bad poles).
For very small~$n$, this will be a rational interpolant,
of degree $\lceil (n-1)/2\rceil$, with error exactly zero on the grid in principle
though nonzero in floating-point arithmetic.  In the figure, the last
such $n$ is marked by a thicker dot.  For most $n$, AAA terminates with
a rational approximant of degree less than
$\lceil (n-1)/2\rceil$ that matches the data to accuracy $10^{-13}$
on the grid without interpolating exactly.
We think of this as a numerical interpolant, since the
error on the grid is so small, whereas much larger errors
are possible between the grid points.  As the grid gets finer, the 
errors between grid points reduce until the tolerance of $10^{-13}$
is reached all across $[-1,1]$.

Another observation about Figure~\ref{figbig} is that
the Floater--Hormann \verb|'equi'| meth\-od is very good~\cite{bos}.
Unlike AAA in its pure form, without partial fractions correction, it
is guaranteed always to produce an interpolant that is pole-free
in $[-1,1]$.

The slowest method to converge is often cubic splines, whose behavior is rock
solid algebraic at the fixed rate $\|f-r\| = O(n^{-4})$, assuming
$f$ is smooth enough.  The convergence of spline
approximations could be speeded up by using degrees that increase
with $n$ (no doubt at the price of some of that rock solidity).

In panels (A) and (D) of Figure~\ref{figbig}, the polynomial least-squares 
approximations converge at first but eventually
diverge exponentially because of unstable amplification of
rounding errors.  Note that the upward-sloping red curves in these two
figures both extrapolate back to about $10^{-16}$, machine precision; the
dotted red lines mark the prediction $10^{-16}\times (1.14)^n$ from (\ref{kuijlaars}) with $\gamma = 2$.
Before this point, it is interesting to compare the very different initial phases 
for $f_A^{}$, with singularities near
$x=\pm 1$, and $f_B^{}$, with singularities near $x=0$.  Clearly we
have initial convergence in the first case and initial divergence in the second,
a consequence of Runge's principle that convergence of polynomial interpolants
depends on analyticity near the middle of the interval.
The figure for the function $f_E^{}$ looks much like that for $f_B^{}$.

The Fourier extension method, as a rule, does somewhat better than
polynomial least-squares in Figure~\ref{figbig}; in certain
limits one expects Fourier methods to have an advantage over
polynomials of a factor of $\pi/2$~\cite[chapter~22]{ATAP}.  Perhaps not too much should
be read into the precise positions of these curves in
the figure, however, as both methods have been implemented
with arbitrary choices of parameters that might have been adjusted
in various ways.

\section{Convergence properties}

What can be said in general about AAA approximation of equispaced
data?  We shall organize the discussion around two questions to be
taken up in successive subsections.
\vspace{7pt}
\begin{itemize}
\item
How does the method normally behave?
\item
How is this behavior consistent with the impossibility theorem?
\end{itemize}
\vspace{7pt}
\noindent
It would be good to support our observations with theorems guaranteeing the success of
the method under appropriate hypotheses,
but unfortunately, like most methods of rational approximation, AAA lacks a theoretical
foundation.

A key property affecting all of the discussion is that, unlike four
of the other five methods of Figure~\ref{figbig} (all but
Floater--Hormann \verb|'equi'|), AAA approximation is
nonlinear.  As so often happens in computational mathematics, the nonlinearity
is essential to its power, while at the same time leading to analytical challenges.  
For example, it means that the theory of frames in numerical approximation,
as presented in~\cite{fna1,fna2,ds}, is not directly applicable.
A theme of that theory, however, remains important here, which is to 
distinguish between {\em approximation\/} and
{\em sampling}.  The approximation issue is, how well can rational
functions approximate a smooth function $f$ on $[-1,1]$?  The sampling
issue is, how effectively will an algorithm based on equally spaced samples find these
good rational approximations?

Our discussion will make reference to the five example functions
$f_A^{},\dots, f_E^{}$ of Figure~\ref{figbig}, which are illustrative of many more 
experiments we have carried out, and in addition we will consider
a sixth function.  In any analysis of polynomial approximations on
$[-1,1]$, and also in the proof of the impossibility theorem
(even though this result is not restricted to polynomial approximations),
one encounters functions analytic inside a {\em Bernstein ellipse} in
the complex plane, which means an ellipse with foci $\pm 1$.  
We define the {\em amber function} (Bernstein means amber
in German) by its Chebyshev series
\begin{equation}
A(x) = \sum_{k=0}^\infty  2^{-k} s_k T_k(x),
\label{amber}
\end{equation}
where the numbers $s_k = \pm 1$ are determined by the binary
expansion of $\pi$,
\begin{equation}
\pi = 11.00100100001111110110\dots\kern 1pt{}_2^{},
\label{pi}
\end{equation}
with $s_k=1$ when the bit is $1$
and $s_k=-1$ when it is $0$.  In Chebfun, one can construct $A$ with the commands
\vspace{6pt}
{\small
\begin{verbatim}
              s = dec2bin(floor(2^52*pi));
              c = 2.^(0:-1:-53)';
              ii = find(s=='0'); c(ii) = -c(ii);
              A = chebfun(c,'coeffs');
\end{verbatim}
\par}
\vspace{6pt}

\begin{figure}
\begin{center}
\kern 7pt
\includegraphics[scale=.7]{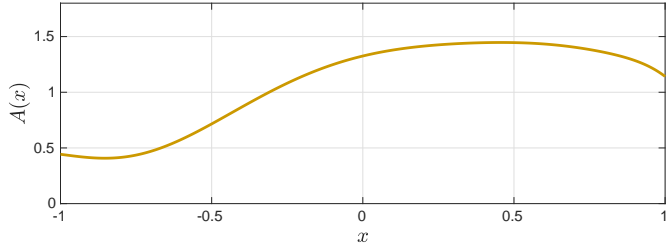}
\caption{\label{figamber}The amber function, a test function constructed to be
analytic in the Bernstein
$2$-ellipse but not in any larger region of
the complex $x$-plane.}
\end{center}
\end{figure}

\noindent
The point of $A(x)$ is that it is analytic in the Bernstein
$2$-ellipse but has
no further analytic structure beyond that, since the bits of $\pi$ are
effectively random.  In particular, it is not analytic or meromorphic in
any larger region of the $x$-plane.  (We believe it has
the $2$-ellipse as a natural boundary~\cite[chap.~4]{kahane}.)
Figure~\ref{figamber} sketches $A(x)$ over $[-1,1]$.

\begin{figure}
\vspace*{10pt}
\begin{center}
\includegraphics[scale=.9]{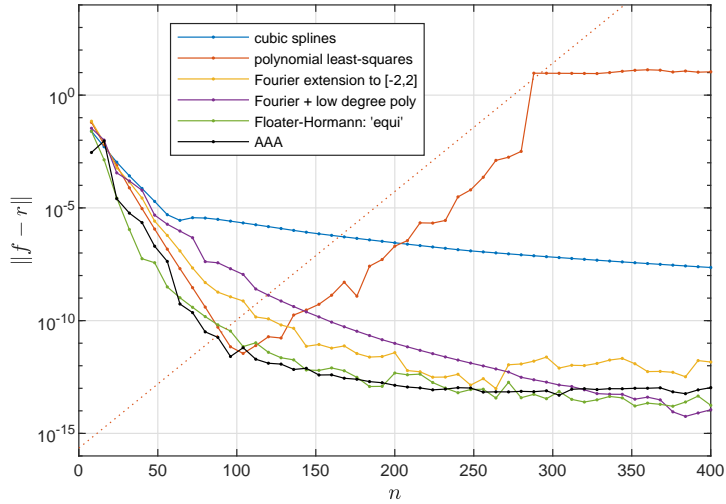}
\caption{\label{ambercurves}A convergence plot as
in Figure~$\ref{figbig}$ for the amber function $A(x)$ of~$(\ref{amber})$.  Note
that for this function, which has no analytic structure beyond analyticity in the
Bernstein $2$-ellipse, the AAA and {\tt 'equi'} methods perform similarly.}
\end{center}
\end{figure}

Figure~\ref{ambercurves} is another plot along the lines of Figure~\ref{figbig}, but
for $A(x)$, and extending to $n=400$ instead of $200$.
We are now prepared to examine the properties of AAA interpolation.

\bigskip
\noindent{\bf 4.1. How does the method normally behave?}
\smallskip

We believe the usual behavior of AAA equispaced interpolation is as follows.  For
small values of $n$, there is a good chance that $f$ will be poorly resolved on the
grid, and the initial AAA interpolant will have poles between the grid points in $[-1,1]$.  In such cases,
as described in the last section, the method switches to a least-squares fit that often
produces acceptable accuracy but without outperforming other methods.  

As $n$ increases, however, $f$ begins to be resolved, and here rational approximation shows its
power.  If $f$ happens to be itself rational, like the Runge function $1/(1+25x^2)$ 
used for experiments in a number of other papers,
AAA may capture it exactly.  More typically, $f$ is 
not rational but, as in the examples of Figure~\ref{figbig}, it has analytic
structure that rational approximants can exploit.  If it is meromorphic, like $\tanh(5x)$,
then AAA quickly finds nearby poles and therefore converges at an accelerating rate.
Even if it has branch point singularities, rapid convergence still takes
place~\cite{clustering}.

\begin{figure}
\begin{center}
\vspace*{8pt}
\includegraphics[scale=.85]{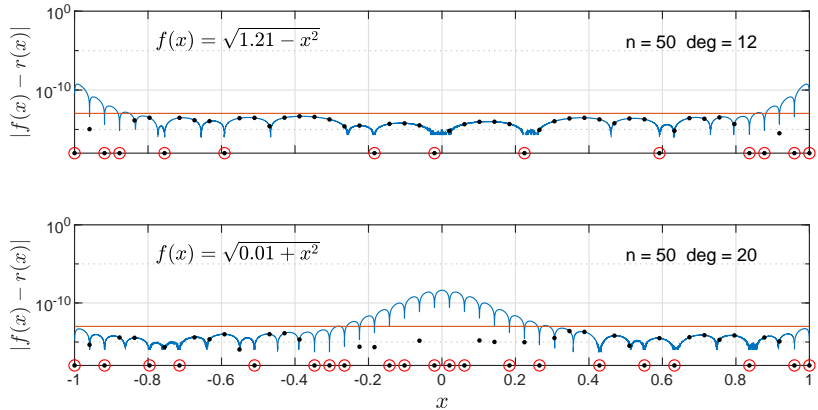}
\caption{\label{twoposs}AAA approximation in its mid-phase of
rapid convergence at $n=50$ for two different functions $f(x)$.
Black dots mark the $n$ sample points, with errors below the AAA relative 
tolerance level of $10^{-13}$ marked by the red line.  The circled black dots
are the subset of AAA support points, where the error in principle is $0$ (apart from
rounding errors), though it has been artificially plotted at $10^{-18}$.
With $f(x) = \sqrt{1.21-x^2}$, above, the big errors between sample points are
near the boundaries, whereas with
$f(x) = \sqrt{0.01+x^2}$, below, they are in the interior.}
\end{center}
\vspace*{-3pt}
\end{figure}

In this middle phase of rapid convergence of AAA approximation, the errors are
many orders of magnitude bigger between the grid points (e.g., $10^{-6}$) than
at the grid points ($10^{-13}$).  The big errors may be near the
endpoints, the pattern familiar in polynomial interpolation since Runge,
but they may also
be in the interior, as happens for example with approximation
of $f(x) = \sqrt{0.01 + x^2}$. Figure~\ref{twoposs} illustrates
these two possibilities.  Convergence eventually happens because the grid points
get closer together and the big errors between them are clamped down.

\begin{figure}
\begin{center}
\kern 8pt
\includegraphics[scale=1.05]{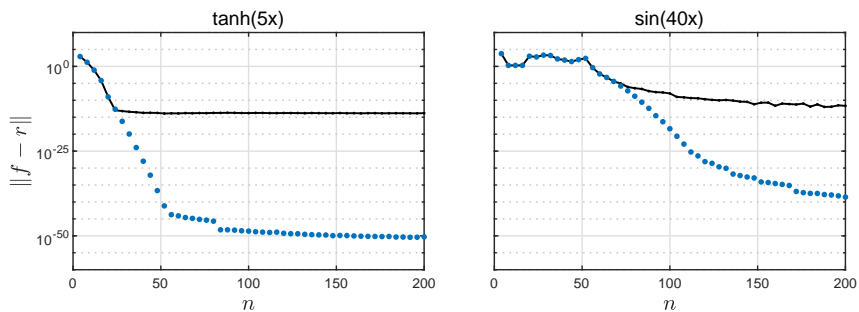}
\caption{\label{julia}AAA errors for two of the functions
of Figure~$\ref{figbig}$ computed in $77$-digit precision
with Julia.  The solid lines are based on the usual AAA tolerance
of $10^{-13}$, and the blue dots are based on
tolerance $10^{-50}$.  
(This computation does not check for poles in $[-1,1]$, which should
in principle lead to values $\|f-r\|=\infty$ especially
in the early stages of the curves on the right; but here the error
is just measured on a $1000$-point equispaced grid.)}
\end{center}
\end{figure}

The AAA method does not keep converging for $n\to\infty$, however.
Instead, it eventually slows down and is limited by its prescribed
relative accuracy, $10^{-13}$ by default.  Thus although it gets high accuracy faster
than the other methods, in the end it too levels off.
To illustrate the significance of the AAA tolerance, Figure~\ref{julia}
repeats the error plots for $f_C^{}(x) = \tanh(5x)$ and $f_D^{}(x)
= \sin(40 x)$ for $n = 4,8,\dots, 200$, but now calculated in 77-digit
BigFloat arithmetic using Julia (with the GenericLinearAlgebra
package) instead of the usual 16-digit
floating point arithmetic.  The solid curve shows behavior with tolerance
$10^{-13}$ and the blue dots with tolerance $10^{-50}$.  

The amber function $A(x)$ was constructed to have no hidden analytic
structure to be exploited; we think of it as being as far from rational as possible.
In Figure~\ref{ambercurves}, this is
reflected in the fact that AAA and polynomial approximants converge at approximately the
same rate until the latter begins to diverge exponentially.
Note also that in Figure~\ref{ambercurves}, unlike the five plots of Figure~\ref{figbig},
AAA fails to outperform the Floater--Hormann \verb|'equi'| method.
This is consistent with the view
that AAA is a robust strategy that exploits analytic structure whereas
Floater--Hormann is a robust interpolation strategy that does not exploit analytic structure.

\bigskip
\noindent{\bf 4.2. How is this consistent with the impossibility theorem?}
\smallskip

In the introduction we summarized the
impossibility theorem of~\cite{ptk} as follows.
In approximation of analytic
functions from $n$ equispaced samples,
exponential convergence as $n\to\infty$ is only possible in
tandem with exponential instability; conversely,
a stable algorithm can converge at best root-exponentially.
The essential reason for this (and the essential construction in
the proof of the theorem) can be summarized in a sentence.
{\em Some analytic functions
are much bigger between the sample points than they are at the sample points;
thus high accuracy requires some approximations to be huge.}
We now explain how the
theorem relates to the six numerical methods presented in Figures~\ref{figbig}
and~\ref{ambercurves}.

{\em Fourier series plus polynomials}, with our choice of polynomial degree
$O(\sqrt n\kern 1pt)$, converge at a root-exponential rate.  It it neither
exponentially accurate nor exponentially unstable.

The {\em Floater--Hormann \verb|'equi'| interpolant} also converges (it appears)
at a root exponential rate, for reasons related to its adaptive choice of degree.

{\em Cubic splines} converge at a lower rate, just $O(n^{-4})$ for
smooth $f$.  Again this method is neither
exponentially accurate nor exponentially unstable.

{\em Fourier extension} also appears to converge root-exponentially, making it, too,
neither exponentially accurate nor exponentially unstable.

{\em Polynomial least-squares,} however, reveals the hugeness of
certain functions.  In the terms of the
theorem, this is the only one of our methods that
appears to be exponentially accurate and exponentially unstable.

\begin{figure}
\begin{center}
\kern 8pt
\includegraphics[scale=1.05]{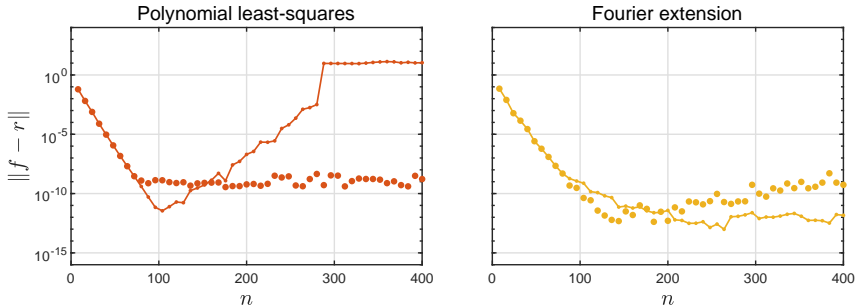}
\caption{\label{twobytwo}The solid lines repeat the polynomial least-squares
(left) and Fourier extension (right) curves of Figure~$\ref{ambercurves}$.
The dots show corresponding results for alternative implementations of
each method.  On the left, switching to the ill-conditioned monomial basis for polynomial
least-squares cuts off the
exponential acuracy and also the exponential instability.  On the
right, switching to the well-conditioned Vandermonde with Arnoldi
basis for Fourier extension
introduces exponential accuracy and exponential instability.}
\end{center}
\end{figure}

Our statements about these last two methods, however, come with a big qualification, which
is illustrated in Figure~\ref{twobytwo}.
In fact, the difference
between Fourier extension and polynomial least-squares lies not in their
essence but in the fashion in which they are implemented.  If you implement
either one with a well-conditioned basis, then it is exponentially accurate
and exponentially unstable.  This is what we have done with the polynomial least-squares
method, which uses the well-conditioned basis of Chebyshev polynomials.  The
Fourier extension method, on the other hand, was implemented with
the ill-conditioned basis of
complex exponentials $\exp(i\pi k x/2)$.  In an ill-conditioned basis like this, high
accuracy will require huge coefficient vectors~\cite{fna1,fna2,ds}, but rounding errors prevent
their computation in floating point arithmetic through the mechanism of
matrices whose condition numbers are unable to grow bigger than
$O((\varepsilon_{\hbox{\scriptsize\rm machine}}^{})^{-1})$.
It is these rounding errors that make our
implementation of Fourier extension stable.  Re-implemented
via Vandermonde with Arnoldi~\cite{vander}, as shown in Figure~\ref{twobytwo},
it becomes exponentially accurate and exponentially unstable.
Conversely, we can implement polynomial least-squares with the exponentially
ill-conditioned monomial basis instead of Chebyshev polynomials.
Because of rounding errors, it then loses its its exponential
accuracy and its exponential instability, as also shown in Figure~\ref{twobytwo}.

Finally, what about AAA?
The experiments suggest it is neither exponentially
accurate nor exponentially unstable.
Insofar as the impossibility theorem is concerned,
there is no inconsistency. Still, what is the mechanism?
In Figure 3.1 we have highlighted the last value of $n$
for which AAA interpolates the data.
In these and many other examples,
that value is very small. Afterwards,
AAA favors closer fits to the data over increasing degrees
of rational approximation, and this has the familiar effect 
of oversampling. Yet, owing to its nonlinear nature,
AAA is free to vary the oversampling factor---and convergence
rates along with it---depending on the data and on the chosen tolerance.

The stability of AAA also stems from the representation
of rational functions in barycentric form, which is discussed
in the original AAA paper~\cite{aaa}.

\section{Discussion}

\begin{figure}
\begin{center}
\vspace*{28pt}
\includegraphics[scale=.7]{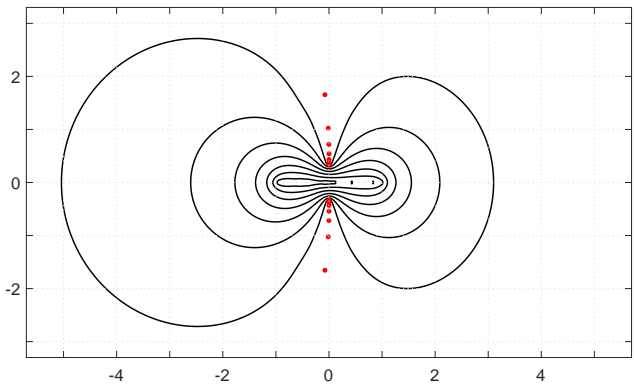}
\caption{\label{fig5}The AAA approximation of Figure~$\ref{fig1}$
extended into the complex plane.  The dots are the poles of the rational
approximant, and the curves are level curves of the error $|f(z)-r(z)|$
with levels (from outside in) $10^{-2}, 10^{-4},\dots, 10^{-14}$.}
\end{center}
\end{figure}

Although we have emphasized just the behavior on $[-1,1]$,
it is well known that rational approximants have good properties of
analytic continuation, beyond the original approximation set.
The AAA method certainly partakes of this advantageous behavior.  
For example, Figure~\ref{fig5} shows the approximation of Figure~\ref{fig1} again
($f(x) = e^x/\sqrt{1+9x^2}$ sampled at $50$ equispaced points in $[-1,1]$), but now evaluated
in the complex plane.  There are many digits of accuracy far from the
approximation domain $[-1,1]$.  This is numerical analytic continuation, and the
other methods we have compared against have no such capabilities.

Another impressive feature of rational approximation is its ability to
handle sampling grids with missing data without much loss of accuracy.
A striking illustration of this effect is presented in Figure~2 of~\cite{wdt}.

The AAA algorithm, as implemented in Chebfun, has a number of adjustable
features.  We have bypassed all these, avoiding both the ``cleanup''
procedure and the Lawson iteration~\cite{aaa} (which is not normally invoked
in any case, by default).

One of the drawbacks of AAA approximation is that although it is 
extremely fast at lower degrees, say, $d< 100$, it slows down
for higher degrees: the complexity is $O(m d^{\kern .8pt 3})$, where $m$ is
the size of the sample set and $d$ is the degree.  For most applications,
we have in mind a regime of problems with $d<100$.
None of the examples shown in this paper come close to this limit.
(With the current Chebfun code, one could write for example
\verb|aaa(F,X,'mmax',200,'cleanup','off','lawson',0)|.)

Our discussion has assumed that the data $f(x_k)$ are accurate samples
of a smooth function, the only errors being rounding errors down at the relative level of
machine precision, around $10^{-16}$.  The Chebfun default tolerance
of $10^{-13}$ was set with this error level in mind.
To handle data contaminated by noise at a higher level $\varepsilon$, we
recommend running AAA with its tolerance parameter \verb|'tol'| set to one or two orders
of magnitude greater than $\varepsilon$.
For unknown noise levels, it should be possible to devise adaptive methods based
on apparent convergence rates---detecting the bend in an L-shaped curve---but we
have not pursued this.  Another approach to dealing with noise in rational approximation
is to combine AAA fitting
in space with calculations related to Prony's method in Fourier space, as advocated
by Wilber, Damle, and Townsend~\cite{wdt}.

\begin{figure}
\begin{center}
\kern 8pt
\includegraphics[scale=.9]{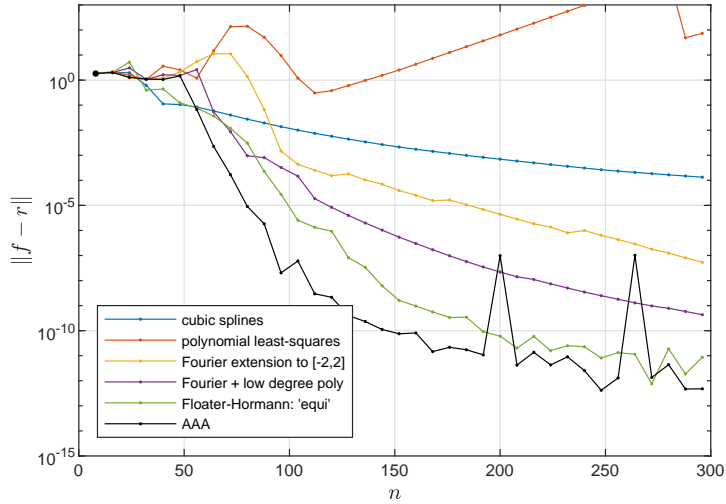}
\caption{\label{all6}AAA approximation of the function $f(x)$ defined
as the sum of all six functions 
$f_A^{},$ $f_B^{},$ $f_C^{},$ $f_D^{},$ $f_E^{},$ and $A$ considered in
this paper.  As usual, AAA mostly outperforms the other methods, but there
are blips at $n=200$ and $n=264$ corresponding to poles in the
approximation interval $[-1,1]$.}
\end{center}
\end{figure}

Most of the methods we have discussed are linear, but AAA is not.  This
raises the question, will it do as well for a truly complicated ``arbitrary'' function
as it has done for the functions with relatively simple properties we have
examined?  As a check of its performance in such a setting, Figure~\ref{all6} repeats
Figure~\ref{ambercurves}, but now for the function $f$ consisting of the
sum of all six test functions we have considered: 
$f_A^{},$ $f_B^{},$ $f_C^{},$ $f_D^{},$ $f_E^{},$ and $A$.
As usual, AAA outperforms the other methods, but blips in
the convergence curve  at $n=200$ and $264$ highlight that it comes with no guarantees.
Both blips correspond to cases where ``bad poles'' have turned up in
the approximation interval $[-1,1]$.  
These problems are related to rounding errors, as can be confirmed
by an implementation in extended precision arithmetic as in Figure~\ref{julia}
or by simply raising the AAA convergence tolerance to $10^{-11}$.
It does seem that further research about avoiding unwanted poles in
AAA approximation is called for, but fortunately, in a
practical setting, such poles are
immediately detectable and thus pose no risk of inaccuracy
without warning to the user.

\section*{Acknowledgments}
We are grateful for helpful suggestions from
John Boyd, Bengt Fornberg, Karl Meerbergen, Yuji Nakatsukasa,
and Olivier S\`ete.

\section*{Appendix: MATLAB/Chebfun code for Figure 2}

\indent~~
\vspace{10pt}

{\footnotesize
\begin{verbatim}
xx = linspace(-1,1,1000)';
MS = 'markersize'; FS = 'fontsize'; LW = 'linewidth';
IN = 'interpreter'; LT = 'latex'; warning off
for plt = 1:5
   switch plt
      case 1, s = '(A)  sqrt(1.21-x^2)';
      case 2, s = '(B)  sqrt(0.01+x^2)';
      case 3, s = '(C)  tanh(5*x)';
      case 4, s = '(D)  sin(40*x)';
      case 5, s = '(E)  exp(-1/x^2)';
   end
   f = inline(vectorize(s(6:end)));
   jplot = 2 - mod(plt,2); iplot = ceil(plt/2);
   axes('position',[-.25+.35*jplot 1.04-.32*iplot .3 .27])
   nmax = 200; nn = 4:4:nmax;

% Cubic splines
   nvec = []; errvec = [];
   for n = nn
      X = linspace(-1,1,n)';
      p = chebfun.spline(X,f(X));
      err = norm(f(xx)-p(xx),inf); nvec = [nvec n]; errvec = [errvec err];
   end
   semilogy(nvec,errvec,'.-',LW,.5,MS,4), grid on, hold on

% Polynomial least-squares, oversampling ratio 2
   nvec = []; errvec = [];
   for n = nn
      X = linspace(-1,1,n)';
      A0 = chebpoly(0:n/2); A = A0(X);
      c = A\f(X); AA = A0(xx); pp = AA*c;
      err = norm(f(xx)-pp,inf); nvec = [nvec n]; errvec = [errvec err];
   end
   semilogy(nvec,errvec,'.-',LW,.5,MS,4), grid on

% Fourier extension on [-2,2], oversampling ratio 2
   nvec = []; errvec = [];
   zz = exp(.5i*pi*xx);
   for n = nn
      X = linspace(-1,1,n)'; Z = exp(.5i*pi*X);
      d = ceil(n/4); c = [real(Z.^(0:d)) imag(Z.^(1:d))]\f(X);
      c = c(1:d+1) - 1i*[0; c(d+2:2*d+1)];
      pp = real((zz.^(0:length(c)-1))*c);
      err = norm(f(xx)-pp,inf); nvec = [nvec n]; errvec = [errvec err];
   end
   semilogy(nvec,errvec,'.-',LW,.5,MS,4), grid on, hold on

% Fourier plus polynomial of degree sqrt(n), oversampling ratio 2
   nvec = []; errvec = [];
   for n = nn
      degp = round(sqrt(n)) - 1;     % degree of polynomial term
      degp = degp + mod(degp+n+1,2); %   ... with opposite parity to n
      degf = (n-1-degp)/2;
      degf = ceil(degf/2);           % oversampling ratio about 2
      X = linspace(-1,1,n)'; 
      C = chebpoly(0:degp);
      A = [cos(pi*(1:degf).*X) sin(pi*(1:degf).*X) C(X)]; c = A\f(X);
      p = @(x) reshape([cos(pi*(1:degf).*x(:)) ...
                        sin(pi*(1:degf).*x(:)) C(x(:))]*c,size(x));
      err = norm(f(xx)-p(xx),inf); nvec = [nvec n]; errvec = [errvec err];
   end
   semilogy(nvec,errvec,'.-',LW,.5,MS,4), grid on, hold on

% Floater-Hormann rational interpolation: chebfun 'equi'
   nvec = []; errvec = [];
   for n = nn
      X = linspace(-1,1,n)';
      p = chebfun(f(X),'equi');
      err = norm(f(xx)-p(xx),inf); nvec = [nvec n]; errvec = [errvec err];
   end
   semilogy(nvec,errvec,'.-',LW,.5,MS,4), grid on

% AAA
   nvec = []; errvec = []; lastinterp = 0;
   for n = nn
      X = linspace(-1,1,n)';
      [r,pol] = aaa(f(X),X,'cleanup','off');
      if any(imag(pol)==0 & abs(pol)<1)           % if there are poles in [-1,1],
         pol = pol(imag(pol)~=0 | abs(pol)>1);    % switch to AAA-least squares
         pX = pol.'; for j = 1:length(pol); pX(j) = min(abs(pol(j)-X)); end
         A = [X.^0 pX./(X-pol.')]; c = A\f(X);
         r = @(x) reshape([x(:).^0 pX./(x(:)-pol.')]*c, size(x));
      end
      if length(pol) == ceil((n-1)/2), lastinterp = lastinterp+1; end
      err = norm(f(xx)-r(xx),inf);
      nvec = [nvec n]; errvec = [errvec err];
   end
   semilogy(nvec,errvec,'.-k',LW,.5,MS,4), grid on
   if lastinterp > 0, plot(nvec(lastinterp),errvec(lastinterp),'.k',MS,8), end

   set(gca,FS,5), axis([0 nmax 1e-16 1e3])
   ylabel('$\|f-r\|$',FS,7,IN,LT), xlabel('$n$',FS,7,IN,LT)
   if plt < 4, set(gca,'xticklabel',{}), xlabel(' '), end
   if mod(plt,2) == 0, set(gca,'yticklabel',{}), ylabel(' '), end
   if plt == 1 | plt == 4
      plot([0 n],eps*1.14.^[0 n],':','color',[.85 .325 .098])
   end
   s2 = strrep(s,'*','');
   if plt == 4, text(121,10,s2,FS,6), else text(25,10,s2,FS,6), end
   set(gca,'ytick',10.^(-15:5:0),'yminorgrid','off'), hold off, shg
end
axes('position',[.49 .14 .25 .15])
for j = 1:5, plot([0 1],[0 0],'.-'), hold on, end
plot([0 1],[0 0],'.-k'), axis([2 3 2 3]), axis off
legend('cubic splines','polynomial least-squares','Fourier extension to [-2,2]',...
    'Fourier + low degree poly','Floater-Hormann: ''equi''','AAA','location',...
    'southeast',FS,5)
print -depsc testsplt
\end{verbatim}
\par
}

\end{document}